\documentclass[12pt]{article}
\usepackage{amssymb}
\usepackage{epsfig}
\usepackage{amsfonts}
\usepackage{verbatim}
\usepackage{latexsym}
\usepackage{graphicx,amsmath}
\usepackage{color}
\usepackage{amsfonts, amssymb, epsfig,subfigure}
\usepackage{mathrsfs}
\usepackage{stmaryrd}
\usepackage{graphicx,amsmath}
\usepackage{latexsym}
\usepackage{verbatim}
\usepackage{color}
\newtheorem{theorem}{Theorem}[section]

\newtheorem{assp}{Assumption}[section]
\newtheorem{lemma}[theorem]{Lemma}
\newtheorem{rem}{Remark}[section]

\newtheorem{expl}[theorem]{Example}

\makeatletter
\@addtoreset{equation}{section}

\newcommand{\beq}[1]{\begin{equation} \label{#1}}
\newcommand{\eeq}{\end{equation}}
\newcommand{\bed}{\begin{displaymath}}
\newcommand{\eed}{\end{displaymath}}
\newcommand{\bea}{\bed\begin{array}{rl}}
\newcommand{\eea}{\end{array}\eed}

\newcommand{\barray}{\begin{array}{ll}}
\newcommand{\earray}{\end{array}}

\newcommand{\Spvek}[2][r]{%
  \gdef\@VORNE{1}
  \left(\hskip-\arraycolsep%
    \begin{array}{#1}\vekSp@lten{#2}\end{array}%
  \hskip-\arraycolsep\right)}

\def\vekSp@lten#1{\xvekSp@lten#1;vekL@stLine;}
\def\vekL@stLine{vekL@stLine}
\def\xvekSp@lten#1;{\def\temp{#1}%
  \ifx\temp\vekL@stLine
  \else
    \ifnum\@VORNE=1\gdef\@VORNE{0}
    \else\@arraycr\fi%
    #1%
    \expandafter\xvekSp@lten
  \fi}

\def\sqr#1#2{{\vcenter{\vbox{\hrule height.#2pt
\hbox{\vrule width.#2pt height#1pt \kern#1pt
\vrule width.#2pt} \hrule height.#2pt}}}}

\newcommand{\E}{\mathbb{E}}
\newcommand{\PP}{\mathbb{P}}
\newcommand{\RR}{\mathbb{R}}

\def\SS{{\mathbb{S}}}

\def\tr{\triangle}

\def\F{{\cal F}}

\def\<{\langle} \def\>{\rangle}

\def\1{\oslash} \def\2{\oplus} \def\3{\otimes} \def\4{\ominus}
\def\5{\circ} \def\6{\odot} \def\7{\backslash} \def\8{\infty}
\def\9{\bigcap} \def\0{\bigcup} \def\+{\pm} \def\-{\mp}
\def\[{\langle} \def\]{\rangle}

\def\spec{\hbox{\rm spec}}

\def\diag{\hbox{\rm diag}}

\newcommand{\dis}{\displaystyle}

\def\nn{\nonumber}

\def\bc{\begin{center}}       \def\ec{\end{center}}
\def\ba{\begin{array}}        \def\ea{\end{array}}
\def\be{\begin{equation}}     \def\ee{\end{equation}}
\def\bea{\begin{eqnarray}}    \def\eea{\end{eqnarray}}
\def\beaa{\begin{eqnarray*}}  \def\eeaa{\end{eqnarray*}}
\def\la{\label}

\begin{document}


\title{The Numerical Invariant Measure of Stochastic Differential Equations With Markovian Switching\thanks{Funding: The research of the first author was supported by the National Natural Science Foundation of China (11171056, 11471071, 11671072), the Natural Science Foundation of Jilin Province (20170101044JC), and the Education Department of Jilin Province (JJKH20170904KJ).}}
\author{Xiaoyue Li\thanks{School of Mathematics and Statistics, Northeast Normal University, 5268 Renmin St.,
Changchun, Jilin, 130024, China (lixy209@nenu.edu.cn).}
\and Qianlin Ma\thanks{School of Mathematics and Statistics, Northeast Normal University, 5268 Renmin St.,
Changchun, Jilin, 130024, China, and High School Attached to Capital Normal University-North
Daxing Branch School, 10 Xinghe St., Beijing, 102600, China (mglangel@126.com).}
\and Hongfu Yang\thanks{School of Mathematics and Statistics, Northeast Normal University, 5268 Renmin St.,
Changchun, Jilin, 130024, China (yanghf783@nenu.edu.cn).}
\and Chenggui Yuan\thanks{Department of Mathematics, Swansea University, Swansea, Wales SA2 8PP, UK (C.Yuan@swansea.ac.uk).}}
 \date{}
\maketitle

\begin{abstract}
The existence and uniqueness of the numerical invariant measure of  the backward Euler-Maruyama method for stochastic differential
equations  with Markovian switching is yielded, and  it is  revealed that the numerical invariant measure converges to the underlying invariant measure   in the Wasserstein metric. The global Lipschitz condition on the drift coefficients required by
 [J. Bao, J. Shao, and C. Yuan, Potential Anal., 44 (2016), pp. 707-727] and [X. Mao, C. Yuan, and G. Yin, J. Comput. Appl. Math., 174 (2005), pp. 1-27]  is released.  Under a polynomial growth condition imposed on  drift coefficients we show that the convergence is    exponential.  Several examples and numerical experiments are given to verify our theory.

\end{abstract}

\section{\large\bf Introduction }
As one of the important classes of  hybrid systems, stochastic differential
equations (SDEs) with Markovian switching
 have been widely used in
  biology, control problems, neutral activity, mathematical finance and other sciences (see, e.g., the monographs \cite{Mao2006, Yin} and the references therein).
    So far, various dynamical properties including moment boundedness, stability,  ergodicity, recurrence and  transience on SDEs with Markovian switching have  been investigated extensively, refer to \cite{Bakhtin2012,Bao,Ba, Mao2006,Pinsky1992,Shao2013,Shao2015,Yin2,Yin}. Yin and Zhu  \cite[pp.181-280]{Yin},   and Mao and
Yuan \cite[pp.164-190]{Mao2006} investigated the stability of SDEs with Markovian switching and showed that the Markov chain facilitates the stochastic stabilization
in which the stationary distribution
of the Markov chain plays an important role.
Pinsky and Scheutzow \cite{Pinsky1992} revealed the fact  that the overall system  may not to be positive
recurrence (resp. transience) even though each subsystem is.
So,  the dynamical behaviors of SDEs with Markovian switching are significantly different from those  of SDEs.

However,
solving the SDEs with Markovian switching is still a challenging task that
  requires using numerical methods or approximation techniques, see, e.g.,  the monographs \cite{Kloeden, Mao2006,Mao2008, Yin}.
   Some long-time behaviors of the SDEs with Markovian switching, for instance, the almost sure stability and the moment stability,  have been preserved  by the numerical solutions,
  see,  e.g., \cite{Higham07, Mao2006, Mao3, Pang, Yin, Zong} and the references therein.    For deterministic systems, the stability of  equilibrium point is among of the interesting topics.   However, many stochastic systems don't posses a deterministic  equilibrium
 state.
Recently, for stochastic
systems with Markovian switching, the stability of the ``stochastic  equilibrium state''-the existence of the invariant measure has drawn increasing attention \cite{Bakhtin2012,Bao, Shao2013,Shao2015,Yin2,Yin}.  Since
the corresponding Kolmogorov-Fokker-Planck equations are always computationally intensive, it is important to be able to approximate the invariant measure numerically.
 Therefore, approximations of
invariant measures for SDEs with Markovian switching have attracted much attention recently.  Mao et al. \cite{Mao2005}, Yuan and Mao \cite{Yuan3} and  Bao et al. \cite{Bao} made use of Euler-Maruyama (EM) method with a constant step size to approximate the underlying invariant measure while
Yin and Zhu \cite[p.159-179]{Yin} did that using  the EM scheme with the  decreasing step size. In the mentioned papers, both the drift coefficients and the diffusion coefficients of the SDEs with regime switching are required to be global Lipschitz continuous. Although the classical Euler-Maruyama (EM) method is convenient for computations and implementations,  the absolute moments of its approximation for SDEs with super-linear coefficients may diverge to infinity
at a finite time  (see, e.g. \cite{Hut01}). It is well know, see \cite{Ma02}, that the  EM numerical solutions fail  to be ergodic, even when the underlying SDE is geometrically
ergodic.
Many implicit methods were used to study the numerical solutions to
SDEs with nonlinear coefficients (see, e.g., \cite{H01,Mil05}). Higham et al.  \cite{H01} proved that the implicit EM numerical solutions converge strongly to the exact solutions of SDEs with  globally one-sided Lipschitz continuous drift term  and globally Lipschitz continuous diffusion term, but the explicit EM method fails to do that.
  Mattingly et al. \cite{Ma02} introduced variants of the implicit EM method to preserve the ergodicity for SDEs with additional noises usually established through the use of Foster-Lyapunov conditions in \cite{Me1,Me2, Me3}  while Liu and Mao \cite{Liu15} took advantage of the implicit EM method to approximate the stability in distribution of non-globally Lipschitz continuous SDEs.
For the background on the implicit methods, we refer the reader to the books
\cite{Kloeden, Mil02}. Shardlow and Stuart \cite{Shard+stu} established the perturbation theory of geometrically ergodic Markov chain with an application to numerical approximations.

Motivated by the papers above,
  this paper  focuses on using  the backward Euler-Maruyama (BEM) method to approximate the invariant measure of  nonlinear SDEs with Markovian switching that the drift coefficients need not to satisfy the global Lipschitz condition.
   The BEM scheme, which is  implicit in the drift term,   has been implemented for SDEs with Markovian switching  to
   investigate the strong convergence and the approximation of the almost sure stability as well as the moment  stability (see, e.g., \cite{Mao3, Zong, Zhou2} and the
references therein).
The main aim of this paper is to
  study the existence and uniqueness of the numerical
invariant measure of the BEM method and the convergence   in the Wasserstein metric to the  invariant measure of the corresponding exact solution  as well as the convergence rate.

The rest of our paper is organized as follows. Section \ref{sec:per} gives some preliminary results on the existence and uniqueness of the  invariant measure for the exact solution. Section \ref{S-dist} focuses on the existence and uniqueness of the numerical invariant measures in  BEM scheme.
Then we go further to reveal that the numerical invariant measure  converges in the Wasserstein distance to the underlying one.
 Section \ref{sec:example} presents several examples and numerical experiments to illustrate our results.

\section{Preliminary}\label{sec:per}
Throughout this paper,  let   $|\cdot|$ denote the Euclidean norm in $\RR^n:=\RR^{n\times 1}$ and the trace norm in $\RR^{n\times m}$.
 If $A$ is a vector or matrix, its transpose is denoted by $A^T$ and its trace norm is denoted by $|A|=\sqrt{\mathrm{trace}(A^TA)}$.
 For vectors  or matrixes $A$ and $B$ with compatible dimensions, $A B$ denotes the usual matrix multiplication. We  denote the indicator function of a set $\mathbb{D}$ by $I_{\mathbb{D}}$, and $\mathbf{0}\in \RR^n$ is a zero vector. For any $\xi=(\xi_{1},\xi_{2},\cdots,\xi_{n})^T\in \RR^n$, $\xi\gg\mathbf{0}$ means each component $\xi_{j}>0,  j=1,2,\cdots,n$. Define $ \hat{\xi}=\min_{1\leq j\leq n}\xi_j$ and $ \check{\xi}=\max_{1\leq j\leq n}\xi_j$. For any  $a,b\in\RR$,  $a\vee b=\max\{a,b\}$, and $a\wedge b=\min\{a,b\}$. For each $R>0$, let $B_{R}(0)=\{x\in \RR^n:|x|\leq R\}$. Let  $\mathscr{B}(\RR^n)$ denote the family of all Borel sets in $\RR^n$.

Let $( \Omega,~\cal{F}$,~$\PP )$ be a complete probability space, and $\mathbb{E}$ denotes the expectation corresponding to $\PP$. Let~$B(t)$ be an $m$-dimensional Brownian motion defined on this probability space. Suppose that  $\{r(t)\}_{t\geq0}$ is a right-continuous  Markov chain with finite  state space $\SS=\{1, 2, \cdots, N\}$ and independent of the Brownian motion $B(\cdot)$, where $N$ is a positive integer. Suppose  ${  \{ {\cal{F}}_{t}}\}_{t\geq 0} $ is a filtration defined on this probability space satisfying the usual conditions (i.e., it is right continuous in $t$ and $\mathcal{F}_0$ contains all $\mathbb{P}$-null sets) such that $B(t)$ and $r(t)$ are ${\cal{F}}_t$ adapted. The  generator of $\{r(t)\}_{t\geq0}$  is  denoted by $Q=(q_{lj})_{N\times N}$, so that for a sufficiently small $\epsilon>0$,
\begin{align*}
 \mathbb{P}\{r(t+\epsilon)=j|r(t)=l\}=\left\{
\begin{array}{ll}
q_{lj}\epsilon+o(\epsilon),&~~\mathrm{if}~~l\neq j,\\
1+q_{ll}\epsilon+o(\epsilon),&~~\mathrm{if}~~l=j.
\end{array}
\right.
\end{align*}
Here $q_{lj}\geq 0$ is the transition rate from $l$ to $j$
if $l\neq j$ while
  $
q_{ll}=-\sum_{l\neq j}q_{lj}.
$
It is well known that almost every sample path of $r(t)$ is a right-continuous step function with a finite number of simple jumps in any finite subinterval of $\mathbb{R}_+:=[0, +\infty)$ (see \cite[p.17-18]{Anderson}). As a standing hypothesis, we assume that the transition probability matrix $Q$ are  {\it irreducible and conservative}.  So  Markov chain  $\{r(t)\}_{t\geq0}$  has a unique stationary distribution $\mu:=(\mu_{1},\mu_{2},\cdots,\mu_{N}) \gg\mathbf{0}\in \mathbb{R}^{1\times N}$ which can be determined by solving the linear equation
 \begin{align}\label{eq:a1.2}
\mu Q=0,\ \ \ \ \ \ \mathrm{subject\ to}\ \  \ \sum_{j=1}^{N}\mu_j=1.
\end{align}

In this paper, we consider  the two-component diffusion process $(Y(t), r(t))$ described by the SDE with Markovian switching
\begin{align}\label{1.1}
dY(t)=f(Y(t),  r(t))dt+g(Y(t), r(t))\mathrm{d}B(t)
\end{align}
on $t\geq 0$ with the initial data $(Y(0),r(0))=(x, i)\in \RR^n\times \SS$, where
$$f:\RR^n\times \SS \rightarrow \RR^n~~~\mathrm{and}~~~~g:\RR^n\times \SS \rightarrow \RR^{n\times m}.$$ 
For convenience we further impose the following hypothesises on the drift and diffusion  coefficients.
\begin{assp}\label{a1}
 For any ~$j\in \mathbb{S}$, there exists a constant~$\alpha_{j}\in\RR$ such that
 \be\la{a-e1}
           (u-v)^{T}(f(u,j)-f(v,j))\leq\alpha_{j}|u-v|^{2},~~~~\forall u,~v\in\RR^n.
\ee
Moreover, for any~$R\geq0$, there exists a positive constant~$K_R$ such that
$$
|f(u,j)-f(v,j)|\leq K_R|u-v|,
$$
for any~$u,~v\in\RR^n$,~$|u|\vee|v|\leq R$~,~$j\in \mathbb{S}$.
\end{assp}

\begin{assp}\label{a3}
For any~$j\in \mathbb{S}$, there exist  constants~ $h_j\in \RR$ and $h>0$ such that
\be\la{a-e3}
|u-v|^{2}|g(u,j)-g(v,j)|^{2}-2|(u-v)^{T}(g(u,j)-g(v,j))|^{2}\leq h_{j}|u-v|^{4}, \ee and
\be\la{a-e2}
|g(u,j)-g(v,j)|^{2}\leq h |u-v|^{2},
\ee
for any~$u,~v\in\RR^n$.
\end{assp}

Next, for convenience, define
\begin{align}\la{2.6}
\beta_j=2\alpha_{j}+h_{j},~~\beta=(\beta_1, \cdots, \beta_N)^T,~~\lambda=|  \mu  \beta |. 
\end{align}
 Assumptions \ref{a1}  and the elementary inequality imply  that for any $u\in \RR^n$
\begin{align}\label{b1}
\begin{split}
   u^{T}f(u,j)  & \leq \alpha_{j}|u|^{2}+|u^{T}f(0,j)|
 \leq \alpha_{j}|u|^{2}+ \frac{\lambda|u|^{2}}{8}+ \frac{2|f(0,j)|^{2}}{\lambda}\\
& \leq \Big(\alpha_{j}+\frac{1}{8}\lambda\Big)|u|^{2}+\sigma_{1},
\end{split}
  \end{align}
and Assumption \ref{a3} and the elementary inequality imply that
\begin{align}\label{b3}
  |g(u,j)|^{2}\leq2h|u|^{2}+\sigma_{2},
\end{align}
 where $\sigma_{1}= 2\max\limits_{j\in \SS}\{|f(0,j)|^{2}/\lambda\}$ and  $\sigma_{2}=2\max\limits_{j\in \SS}\{|g(0,j)|^{2}\}$. Moreover, choosing constants $2p\leq \varepsilon=\lambda/ 16h $, we find that
  \begin{align}\label{b5}
 |u|^{2}|g(u,j)|^{2}+(p-2)|u^{T}g(u,j)|^{2}
 \leq&\big(h_{j}+(3\varepsilon+2p)h\big)|u|^{4}+\frac{\sigma_2(1+2p+3\varepsilon^{-1})}{2}|u|^{2}\nn\\
 \leq&\big(h_{j}+\frac{1}{4}\lambda\big)|u|^{4}+\sigma_3|u|^{2},
  \end{align}
where $\sigma_3=(1+\lambda/16h+48h/\lambda)\sigma_2/2$.

Under Assumptions \ref{a1} and \ref{a3},   the equation (\ref{1.1}) admits a unique    solution   $(Y(t),r(t))$
(see, e.g., \cite[Theorem 3.17, p.93]{Mao2006}). Throughout the paper, we write $(Y_t^{x, i},r_t^{i})$ in lieu of
$(Y(t),r(t))$ to highlight the initial data $(Y(0),r(0))=(x,i)$.
  Let $\mathcal P(\RR^n\times \SS)$ denote the family of all probability measures on $\RR^n\times \SS$.
  For any $p\in(0,1]$, define a metric on $\RR^n\times \SS$ as the following
$$
d_{p}((u,j),(v,l)):=|u-v|^{p}+I_{\{j\neq l\}},~~~(u,j),(v,l)\in \RR^n\times \SS,
$$
and the corresponding  Wasserstein distance between $\nu,\tilde{\nu}\in \mathcal P(\RR^n\times \SS)$ by
\begin{align*}
W_{p}(\nu,\tilde{\nu}):=\inf_{\pi\in C(\nu,\tilde{\nu})}\int_{  {(\RR^{n}\times \SS) \times(\RR^{n}\times \SS)} } d_{p}(u,v)\pi(\mathrm{d}u,\mathrm{d}v),
\end{align*}
where $C(\nu,\tilde{\nu})$ denotes the set of all couplings of $\nu$ and $\tilde{\nu}$. Let $\mathbf{P}_{t}(x,i;\mathrm{d}u\times \{l\})$ be the transition probability kenel of the pair $\big(Y_t^{x,i},r^{i}_t\big)$, a time homogeneous Markov process (see, e.g, \cite[Theorem 3.28, pp.105-106]{Mao2006}). Recall that $\pi\in\mathcal P(\RR^n\times S)$ is called an invariant measure of $\big(Y_t^{x,i},r_t^{i}\big)$ if
\begin{align*}
\pi(\Gamma\times\{j\})=\sum_{l=1}^{N}\int_{\RR^n} \mathbf{P}_{t}(u,l;\Gamma\times\{j\})\pi(\mathrm{d}u\times\{l\}),~~~\forall t\geq0,~\Gamma\in\mathscr{B}(\RR^n),~j\in \SS
\end{align*}
holds.  For each~$p>0$, define
\begin{align}\label{4.24}
  \Lambda=\diag(8\beta_1+7\lambda, \cdots, 8\beta_{N}+7\lambda),~~  Q_{p}=Q+\frac{p}{16}\Lambda,~~
\eta_{p}=-\max\limits_{\gamma\in \spec(Q_{p})}\mathrm{Re}\gamma, 
\end{align}
where~$\lambda$ and $\beta_j $ are introduced in \eqref{2.6}, $Q$ is the generator  of $\{r(t)\}_{t\geq 0}$, and $\spec(Q_{p})$   denotes the spectrum of $Q_{p}$.

 The following lemma highlights the relationship between the sign of $\mu \beta$ and the sign of $\eta_{p}$.
\begin{lemma}\la{yhf2.1}
For any $p>0$, there exists a positive constant  $H(p)$ such that for any $t>0$
\begin{align*}
\mathbb{E}\left[\exp\bigg(\frac{p}{16}
\int_{0}^{t}\Big(8\beta(r(s))+7\lambda\Big)
\mathrm{d}s\bigg)\right]\leq H(p)\mathrm{e}^{-\eta_p t}.
\end{align*}
Moreover, if
$
 \mu \beta<0,$ there is a constant $\bar{p}>0$ such that $\eta_{p}>0$ for  $p\in (0, \bar{p})$. Furthermore,
\begin{description}
  \item[(1)] $\bar{p}=+\infty$   if ~$ 8\check{\beta}+ 7\lambda \leq0$;
    \item[(2)]  $\bar{p}\in\left(0,\min\limits_{j\in \SS,~ 8\beta_{j}+7\lambda>0}\Big\{-16 q_{jj}/(8\beta_{j}+7\lambda)\Big\}\right)$  if ~$ 8\check{\beta}+ 7\lambda >0$,
\end{description}
where $\lambda$ and $\beta_{j}$  are introduced in \eqref{2.6}.
\end{lemma} 
\begin{proof}  According to \eqref{eq:a1.2} and $\mu \beta<0$, it is easy to obtain
  $$\sum_{j=1}^{N}\mu_{j}(8\beta_{j} +7 \lambda)=8\mu\beta+7\lambda=\mu\beta<0.$$
Then the desired assertion follows from \cite[Proposition 4.1 and Proposition 4.2]{Ba} directly.
\end{proof}

We have  the following result on the invariant measure for the exact solution.
\begin{theorem}\la{yth3.1}
 Suppose that Assumptions \ref{a1}, \ref{a3}, and $\mu \beta<0$ hold, then
 the solutions of the SDE with
Markovian switching (\ref{1.1}) converge to a unique invariant measure $\pi\in \mathcal{P}(\RR^n\times \SS)$ with some exponential rate $\xi>0$ in the Wasserstein distance.
\end{theorem}
\begin{proof} We shall adopt the approach  of \cite[Theorem 2.3]{Bao} to complete the proof. Let
{\begin{align}\la{yp}
p_0=1\wedge  \bar{p}\wedge \lambda/32h.
  \end{align}}
Thus, for any  $p\in (0,p_0)$,  \eqref{b5} holds, and $\eta_p>0$ follows from Lemma \ref{yhf2.1}.
One observes that
\begin{eqnarray}\label{eeq+1}
 {\cal{L} } \left(  (1+|x|^{2})^\frac{p}{2}\xi_{i}^{(p)}\right)
 \leq   C - {\eta_{p} } \xi_i^{(p)}  (1+|x|^{2})^\frac{p}{2},
\end{eqnarray}
for $p\in (0,p_0)$, where $\xi^{(p)}=(\xi_{1}^{(p)},\cdots,\xi_{N}^{(p)})\gg\mathbf{0}$ is a eigenvector of
$Q_{p}$ corresponding to $-\eta_{p}$, $C$ is  a  positive constant.
Borrowing the proof method  of \cite[Theorem 2.3]{Bao} we can get the result on the existence and uniqueness of the invariant measure but omit the details to avoid duplication. By the similar way to Theorem \ref{yth3.2}, we yield the exponential convergence rate.
\end{proof}

\begin{rem}\la{re3}
{\rm By virtue of Theorem \ref{yth3.1},  the solution $(Y(t),r(t))$ is ergodic and the transition probability of $(Y(t),r(t))$  converges to its invariant measure  with exponential rate in the Wasserstein distance.  Furthermore, due to (\ref{eeq+1})  the Foster-Lyapunov criterion \cite[Theorem 6.1, p.536]{Me3}  implies that $(Y(t),r(t))$
is exponentially ergodic, provided all compact sets are petite for some skeleton chain.
Thus, this pair is strongly mixing since it is positively Harris-recurrent, see details in \cite[p.881]{Ath}. However more conditions should be imposed on the coefficients of the equation in order for all compact sets are petite for some skeleton chain.}
\end{rem}

\section{Numerical Invariant Measure}
\label{S-dist}
 This section is devoted to the existence and uniqueness of the
  numerical invariant measure of the BEM method  and     approximation of  the
  numerical invariant measure to  the underlying one  in the Wasserstein metric.  In order to define the numerical solution, we need to explain how to simulate a discrete Markov chain, which has been formulated in  \cite[Chapter 4, p.111]{Mao2006}. To make the content  self-contained,  we sketch it here.

Given a stepsize $\Delta>0$ and let
$ P(\Delta)=\left( P_{ij}(\Delta)\right)_{N\times N}=\exp(\Delta Q)$.  The discrete Markov chain $\{r_k,\ k=0, 1, \cdots\}$ can be simulated as follows:  let $r(0)=i$ and give a random pseudo number $\varsigma_1$ obeying the uniform $(0,1)$ distribution. Define
 \begin{align*}
r_1=\left\{
\begin{array}{ll}
\dis i_1,~~~~~\mathrm{if}~~i_1\in \SS-\{N\}~~\mathrm{such\ that}~~\sum_{j=1}^{i_1-1}P_{i  j}(\Delta)\leq \varsigma_1< \sum_{j=1}^{i_1}P_{i  j}(\Delta),\nn\\
\dis N,~~~~~\mathrm{if} ~~\sum_{j=1}^{N-1}P_{i  j}(\Delta)\leq \varsigma_1,
\end{array}
\right.
\end{align*}
where $\sum_{j=1}^{N}P_{i  j}(\Delta)=0$ as usual. In other words, the probability of state $s$ being chosen is given by $\PP(r_1=s)=P_{i  s}(\Delta)$. Generally, after the computations of $r_0, r_1, \cdots, r_k$,  give a random pseudo number $\varsigma_{k+1}$ obeying a uniform $(0,1)$ distribution and define $r_{k+1}$ by
 \begin{align*}
r_{k+1}=\left\{
\begin{array}{ll}
i_{k+1},\! &\dis\mathrm{ if }~i_{k+1}\in \SS-\{N\}~\mathrm{ such\ that }\sum_{j=1}^{i_{k+1}-1}P_{r_k j}(\Delta)\leq \varsigma_{k+1}< \sum_{j=1}^{i_{k+1}}P_{r_k j}(\Delta),\nn\\
N,  &\dis\mathrm{if}~\sum_{j=1}^{N-1}P_{r_k j}(\Delta)\leq \varsigma_{k+1}.
\end{array}
\right.
\end{align*}
This procedure can be carried out independently to obtain more trajectories.

We can now  define the
  BEM scheme for the SDEs   with Markovian switching (\ref{1.1}).   Let  $X_0=x$, $r_0=i$, and define
\begin{align}\label{4.2}
X_{k+1}=X_{k}+f(X_{k+1},r_{k})\triangle+g(X_{k},r_{k})\triangle B_{k},~~k\geq 0,
\end{align}
where $\triangle B_{k}=B(t_{k+1})-B(t_{k})$. Here $X_k, r_k, k\ge 0,$  depend on the step size $\tr$, we drop it  for simplicity.
 We point out that the BEM method (\ref{4.2}) is well-defined under Assumption \ref{a1} based on a known result {\rm \cite[Lemma 5.1]{Mao3}} as follows.
\begin{lemma}
 Let Assumption 2.1 holds and $\triangle< 1/|\check{\alpha}| $. Then for any $j\in \SS$,    $b\in \RR^{n}$, there is a unique root $u\in \RR^{n}$ of the equation
$$
u=b+f(u,j)\triangle.
$$
\end{lemma}
It is useful to write (\ref{4.2}) as
\begin{align}\label{c1}
X_{k+1}-f(X_{k+1},r_{k})\triangle=X_{k}+g(X_{k},r_{k})\triangle B_{k}.
\end{align}
For any $j\in \mathbb{S}$, define a function $G_{j}:\RR^{n}\rightarrow \RR^{n}$ satisfying $G_{j}(u)=u-f(u,j)\triangle$. Then $G_j$ has its inverse function $G^{-1}_{j}:\RR^{n}\rightarrow \RR^{n}$ for any $j\in\mathbb{ S}$. Moreover, the BEM method (\ref{4.2}) can be represented as
\begin{align}\label{c2}
X_{k+1}=G^{-1}_{r_{k}}(X_{k}+g(X_{k},r_{k})\triangle B_{k}),~~~~\forall k\geq 0.
\end{align}

Similar to that of \cite[Theorem 6.14, p.250]{Mao2006}, we can prove the following result.
\begin{lemma}
$\{(X_{k},r_{k})\}_{k\geq0}$ is a time homogeneous Markov chain.
\end{lemma}

Let $\mathbf{P}_{k\triangle}^{\Delta}(x,i;\mathrm{d} u\times\{l\})$ be the transition probability kernel of the pair $\big(X_{k}^{x,i},r_{k}^{i}\big)$, a time homogeneous Markov chain. If $\pi^{\triangle}\in\mathcal P(\RR^n\times \mathbb{S})$ satisfies
\begin{align*}
\pi^{\triangle}(\Gamma\times\{j\})
=\sum_{l=1}^{N}\int_{\RR^n}\mathbf{P}_{k\triangle}^{\triangle}(u,l;\Gamma\times\{j\})\pi^{\triangle}(\mathrm{d}u\times\{l\}),
\forall k\geq0, \Gamma\in\mathscr{B}(\RR^n), j\in \SS,
\end{align*}
then $\pi^{\triangle}$ is called an invariant measure of $\big(X_{k}^{x,i},r_{k}^{i}\big)$. For convenience,  Denote by $C$  a generic positive constant which value may be different  with different  appearance and is independent of  the iteration number $ k$  and the time stepsize $\triangle$.

In order to show  the existence of the numerical invariant measure we  prepare the following lemma on the moment boundedness of the numerical solution of the BEM scheme borrowing the idea of \cite{Liu15}.
\begin{lemma}\la{le-mom}
Under the conditions of  Theorem \ref{yth3.1}, there exists a constant  $\bar{\tr} $  such that the numerical solution of BEM scheme with any initial value $(x,i)\in\RR^n\times \SS$ satisfies
\begin{align}\label{4.13}
\sup\limits_{k\geq0}\mathbb E|X_{k}|^{p}
\leq C(1+|x|^p)
\end{align}
for any $\tr\in (0,\bar{\tr} )$ and any $p\in(0,p_0)$, where  $p_0$ is defined by \eqref{yp}. 
\end{lemma}
\begin{proof}
It follows from (\ref{b1}) and (\ref{4.2})  that
\begin{align*}
|X_{k+1}|^{2}
 =& X_{k+1}^{T}\Big(f(X_{k+1},r_{k})\triangle+X_{k}+g(X_{k},r_{k})\triangle B_{k}\Big)\nn\\
 \leq&\big(\alpha_{r_{k}}+\frac{1}{8}\lambda\big)|X_{k+1}|^{2}\triangle+\sigma_{1}\triangle+\frac{1}{2}|X_{k+1}|^{2}+\frac{1}{2}|X_{k}+g(X_{k},r_{k})\triangle B_{k}|^2.
\end{align*}
Choosing a constant  $0<\tr_1<1$ such that $(  2 \breve{| \alpha |} +\frac{1}{4}\lambda) \triangle_1\leq 1/3$ (where $\breve{| \alpha |}:= \min_{i\in \SS}{| \alpha_i |}$),  we then  obtain
for any $\triangle \in (0,  \triangle_1]$,
\begin{align*}
|X_{k+1}|^{2}\leq\frac{1}{1-(2\alpha_{r_{k}}+\frac{1}{4}\lambda )\triangle}|X_{k}+g(X_{k},{r_{k}})\triangle B_{k}|^2+\frac{2\sigma_{1}\triangle}{1-(2\alpha_{r_{k}}+\frac{1}{4}\lambda )\triangle},
\end{align*}
which implies
\begin{align*}
1+|X_{k+1}|^{2}
\leq &\frac{1}{1-(2\alpha_{r_{k}}+\frac{1}{4}\lambda )\triangle}\Big[1+|X_{k}+g(X_{k},{r_{k}})\triangle B_{k}|^2+\Big(2\sigma_{1}-2\alpha_{r_{k}} \Big)\triangle\Big]\nn\\
\leq  &\frac{ (1+|X_{k}|^{2})}{1-(2\alpha_{r_{k}}+\frac{1}{4}\lambda)\triangle}\Big(1+ \upsilon_{k}(r_{k})\Big),
\end{align*}
where
\begin{align*}
\upsilon_{k}(r_{k})=\frac{2X_{k}^{T}g(X_{k},r_{k})\triangle B_{k}+|g(X_{k},r_{k})\triangle B_{k}|^{2}+c_1\triangle}{1+|X_{k}|^{2}},~~~ c_1=|2\sigma_{1}-2\hat{\alpha} | .
\end{align*}
For any $p\in(0, p_0)$ where $p_0$ is defined by \eqref{yp},  noting that
\begin{align}\label{4.5}
(1+u)^{\frac{p}{2}}\leq1+\frac{p}{2}u+\frac{p(p-2)}{8}u^{2}+\frac{p(p-2)(p-4)}{48}u^{3},~~u\geq -1
\end{align}
and $\upsilon_{k}(r_{k})>-1$,
we then have
\begin{equation}\label{4.6}
\begin{split}
  &\mathbb E\Big((1+|X_{k+1}|^{2})^{\frac{p}{2}}|\mathcal F_{t_k}\Big)
\leq  \frac{ (1+|X_{k}|^{2})^{\frac{p}{2}}}{[1-(2\alpha_{r_{k}}+\frac{1}{4}\lambda)\triangle]^{\frac{p}{2}}} \\
&~~~~\times \mathbb E\Big(1+\frac{p}{2}\upsilon_{k}(r_{k})+\frac{p(p-2)}{8}\upsilon_{k}^{2}(r_{k})
 +\frac{p(p-2)(p-4)}{48}\upsilon_{k}^{3}(r_{k})\Big|\mathcal F_{t_k}\Big).
 \end{split}
\end{equation}
Since $\triangle B_{k}$ is independent of $\mathcal F_{t_k}$, we have
$  \mathbb{E }(\triangle B_k|\mathcal F_{t_k}) =0,$ $ \mathbb{E} (|A\triangle B_k|^2|\mathcal F_{t_k})
=|A|^2 \triangle, $
for any $A\in \RR^{n\times m}$.
Hence,
\begin{align} \label{4.7}
 \mathbb E\big(\upsilon_{k}(r_{k})|\mathcal F_{t_k}\big)
 = \frac{|g(X_{k},r_{k})|^{2}\triangle+c_{1}\triangle}{1+|X_{k}|^{2}}.
\end{align}
Using  the properties
$$  \mathbb{E} (|\triangle B_k|^{2j})= C\tr^j,~ \mathbb{E }(|\triangle B_k|^{2j-1}|\mathcal F_{t_k})\leq C\tr^{j-\frac{1}{2}},~ j=2,3,\cdots,
$$
we compute
\begin{equation}\label{4.8}
  \mathbb E\big(\upsilon_{k}^{2}(r_{k})|\mathcal F_{t_k}\big)
 = \frac{1}{(1+|X_{k}|^{2})^{2}}\Big(4|X_{k}^{T}g(X_{k},r_{k})|^{2}\triangle + C\tr^{\frac{3}{2}}
 \Big)
 \geq \frac{4 |X_{k}^{T}g(X_{k},r_{k})|^{2} \triangle}{(1+|X_{k}|^{2})^{2}},
 \end{equation}
and
\begin{align}\label{4.9}
  \begin{split}   \mathbb E\big(\upsilon_{k}^{3}(r_{k})|\mathcal F_{t_k}\big )
&\leq \frac{9}{(1+|X_{k}|^{2})^{3}}\mathbb{E}\bigg[8|X_{k}^{T}g(X_{k},r_{k})\triangle B_{k}|^3+|g(X_{k},r_{k})\triangle B_{k}|^{6}+c_{1}^3\triangle^3 \Big|\mathcal{F}_{t_k}\bigg] \\
 &\leq  C\tr^{\frac{3}{2}}.\end{split}
 \end{align}
 Combining (\ref{4.6})-(\ref{4.9}) and
 using (\ref{b3}), for any $k\geq 0$ we obtain,
\begin{equation}\label{4.10a}
\begin{split}
& \mathbb E\Big((1+|X_{k+1}|^{2})^{\frac{p}{2}}|\mathcal F_{t_k}\Big)\\
 \leq  &\frac{ (1+|X_{k}|^{2})^{\frac{p}{2}}}{[1-(2\alpha_{r_{k}}+\frac{1}{4}\lambda)\triangle]^{\frac{p}{2}}}
 \bigg\{1+\frac{p}{2}\bigg[\frac{|X_{k}|^{2}|g(X_{k},r_{k})|^{2} + (p-2)|X_{k}^{T}g(X_{k},r_{k})|^{2}}{(1+|X_{k}|^{2})^{2}} \tr\\
 &~~~~~~~~~~~~~~~~~~~~~~~~~~~~~~~~~~~~~~~+\frac{\big(2h+c_{1}\big)|X_{k}|^{2}+\sigma_2+c_{1}}{(1+|X_{k}|^{2})^{2}}\triangle\bigg]+ C\tr^{\frac{3}{2}}\bigg\}.
 \end{split}
\end{equation}
This, together with (\ref{b3}) and (\ref{b5}), implies
\begin{equation}\label{4.10}
\begin{split}
& \mathbb E\Big((1+|X_{k+1}|^{2})^{\frac{p}{2}}|\mathcal F_{t_k}\Big)\\
 \leq  &\frac{ (1+|X_{k}|^{2})^{\frac{p}{2}}}{[1-(2\alpha_{r_{k}}+\frac{1}{4}\lambda)\triangle]^{\frac{p}{2}}}
 \bigg\{1+\frac{p}{2}\bigg[\frac{\big( h_{r_{k}}+\frac{1}{4} \lambda\big)|X_{k}|^{4}+ \sigma_3|X_{k}|^{2} }{ (1+|X_{k}|^{2})^{2}} \tr\\
 & ~~~~~~~~~~~~~~~~~~~~~~~~~~~~~~~~~~~~~~~~~+\frac{\big(2h+c_{1}\big)|X_{k}|^{2}+\sigma_2+c_{1}}{(1+|X_{k}|^{2})^{2}}\triangle\bigg]+ C\tr^{\frac{3}{2}}\bigg\}\\
 \leq  &\frac{ (1+|X_{k}|^{2})^{\frac{p}{2}}}{[1-(2\alpha_{r_{k}}+\frac{1}{4}\lambda)\triangle]^{\frac{p}{2}}}
 \bigg[1
 +\frac{p}{2} \big( h_{r_{k}}+ \frac{1}{4}\lambda\big)   \tr + C\tr^{\frac{3}{2}}
 \bigg]+C\triangle.
 \end{split}
\end{equation}
Choosing a constant $0<\tr_2 \leq \tr_1$ sufficiently small such that
  $$C\tr_2^{\frac{1}{2}} \leq {p\lambda }/8,~~\mathrm{and}~~{27 (p+2)\big(2\breve{| \alpha |} + \lambda/4\big)^2}\tr_2\leq 2{\lambda},$$
 this yields that for any $\tr\in (0,\tr_2]$
 \begin{equation}\label{eq01}\frac{p}{2} \big( h_{r_{k}}+ \frac{1}{4}\lambda\big)   \tr + C\tr^{\frac{3}{2}}\leq \frac{p}{2} \big( h_{r_{k}}+ \frac{1}{2}\lambda\big)\tr\end{equation}
 and
\begin{align}\label{4.11}
\begin{split}
 \Big[1-\big(2\alpha_{r_{k}}+\frac{\lambda}{4}\big)\triangle \Big]^{-\frac{p}{2}}
\leq& 1+\frac{p}{2}\big(2\alpha_{r_{k}}+\frac{\lambda}{4}\big)\triangle +\frac{p(p+2)(2\breve{| \alpha |} +\frac{1}{4}\lambda)^2}{8[1-(2\breve{| \alpha |} +\frac{1}{4}\lambda)\tr_2]^{\frac{p}{2}+2}}\tr^2 \\
\leq& 1+\frac{p}{2}\big(2\alpha_{r_{k}}+\frac{5\lambda}{16}\big)\triangle.
 \end{split}
\end{align}
Then  for any $\triangle \in (0, \triangle_2]$, combining \eqref{4.10}-\eqref{4.11} we obtain
\begin{align*}
\mathbb E\big((1+|X_{k+1}|^{2})^{\frac{p}{2}}|\mathcal F_{t_k}\big)
\leq
   (1+|X_{k}|^{2})^{\frac{p}{2}}\Big[1+\frac{p}{2}\big(2\alpha_{r_{k}}+h_{r_{k}}+\frac{13\lambda}{16}\big)\triangle+C\tr^2\Big]+C\triangle.
\end{align*}
Letting $\bar{\tr}$ be a constant such that  $\bar \triangle \in (0, \triangle_2]$, $C\bar{\tr} \leq p\lambda/32$ and $   ( {   {{\breve{|\beta|}}}  } +\frac{7}{8}\lambda ) \bar{\triangle}<1$ (where $ { { \breve{|\beta|}} }=\max_{i\in \SS} |  \beta_i|$),  we arrive at  for
$\triangle \in (0, \bar \triangle]$
\begin{align}\label{eq-tight}
 \mathbb E\big((1+|X_{k+1}|^{2})^{\frac{p}{2}}|\mathcal F_{t_k}\big)
  \leq \Big[1+\frac{p}{2}\big(\beta_{r_{k}}+\frac{7}{8}\lambda\big)\triangle\Big](1+|X_{k}|^{2})^{\frac{p}{2}}+C\triangle,
  \end{align}
where $\beta_i$ is defined as \eqref{2.6} for each $i\in \SS$.
For any $k\geq 1$, we further compute
\begin{align}\la{yhf3.8}
\begin{split}
 \mathbb E\big((1+|X_{k+1}|^{2})^{\frac{p}{2}}|\mathcal F_{t_{k-1}}\big)
\leq&\Big[1+\frac{p}{2}\big(\beta_{r_{k}}+\frac{7}{8}\lambda\big)\triangle\Big] \mathbb{E}((1+|X_{k}|^{2})^{\frac{p}{2}}|\mathcal F_{t_{k-1}})+C\triangle \\
 \leq&\prod_{j=k-1}^{k}\Big[1+\frac{p}{2}\big(\beta_{r_{j}}+\frac{7}{8}\lambda\big)\triangle\Big](1+|X_{k-1}|^{2})^{\frac{p}{2}} \\
 &+C\tr\Big[1+\frac{p}{2}\big(\beta_{r_{k}}+\frac{7}{8}\lambda\big)\triangle\Big]+C\triangle.
 \end{split}
\end{align}
 Repeating (\ref{yhf3.8}) we obtain
\begin{align*}
\mathbb E\big((1+|X_{k+1}|^{2})^{\frac{p}{2}}|\mathcal F_{0}\big)
 \leq&(1+|X_{0}|^{2})^{\frac{p}{2}}\prod_{j=0}^{k}\Big[1+\frac{p}{2}\big(\beta_{r_{j}}+\frac{7}{8}\lambda\big)\triangle\Big]\nn\\
 &+C\tr\sum_{i=1}^{k}\prod_{j=k-i+1}^{k}\Big[1+\frac{p}{2}\big(\beta_{r_{j}}+\frac{7}{8}\lambda\big)\triangle\Big]+C\tr.
\end{align*}
Hence, for any $k\geq0$, by virtue of the homogeneous property of the Markov chain, taking expectations on both sides yields
\begin{align*}
\mathbb E\big((1+|X_{k+1}|^{2})^{\frac{p}{2}}\big)
\leq&(1+|x|^{2})^{\frac{p}{2}}\mathbb E \Bigg[\prod_{j=0}^{k}\Big(1 +\frac{p}{2}\big(\beta_{r_{j}}+\frac{7}{8}\lambda\big)\triangle \Big)\Big] \nn\\
 &+C\tr   \sum_{i=1}^{k}\E \Big[\E\Big(\prod_{j=k-i+1}^{k}\Big(1+\frac{p}{2}\big(\beta_{r_{j}}+\frac{7}{8}\lambda\big)\triangle\Big)|\F_{k-i}\Big)\Big]+C\tr
 \nn\\
 \leq&(1+|x|^{2})^{\frac{p}{2}}\mathbb E \Bigg[\prod_{j=0}^{k}\Big(1 +\frac{p}{2}\big(\beta_{r_{j}}+\frac{7}{8}\lambda\big)\triangle \Big)\Big] \nn\\
 &+C\tr   \sum_{i=1}^{k}  \E\Big[\prod_{j=1}^{i}\Big(1+\frac{p}{2}\big(\beta_{r_{j}}+\frac{7}{8}\lambda\big)\triangle\Big) \Big] +C\tr.
 \end{align*}
Thus, we have
\begin{equation}\label{111}
\begin{split}
   \mathbb  E\big((1+|X_{k+1}|^{2})^{\frac{p}{2}}\big)
\leq &(1+|x|^{2})^{\frac{p}{2}} \mathbb E\Bigg[\exp\Big(\sum_{j=0}^{k}\log\Big(1 +\frac{p}{2}\big(\beta_{r_{j}}+\frac{7}{8}\lambda\big)\triangle\Big) \Big)\Bigg] \\
 &+C\tr   \sum_{i=1}^{k} \E\Bigg[ \exp\Big(\sum_{j=1}^{i}\log\Big(1+\frac{p}{2}\big(\beta_{r_{j}}+\frac{7}{8}\lambda\big)\triangle\Big)\Big) \Bigg]+C\tr.
 \end{split}
\end{equation}
 Then, by  inequality
$\log(1+u)\leq u,~ \forall u>-1,$
we compute
\begin{equation}\label{111}
\begin{split}
   \mathbb  E\Big((1+|X_{k+1}|^{2})^{\frac{p}{2}}\Big)
\leq &(1+|x|^{2})^{\frac{p}{2}} \mathbb E\Bigg[\exp\bigg(\frac{p}{16}\sum_{j=0}^{k}\big(8\beta_{r_{j}}+7\lambda\big)\triangle \bigg)\Bigg] \\
 &+C\tr   \sum_{i=1}^{k} \E\Bigg[ \exp\bigg( \frac{p}{16}\sum_{j=1}^{i}\big(8\beta_{r_{j}}+7\lambda\big)\triangle \bigg) \Bigg]+C\tr.
 \end{split}
\end{equation}
For any $p\in (0,p_0)$,
Lemma \ref{yhf2.1}
implies that $\eta_{p}>0$ and there exists  a positive constant  $H(p)$ such that
\begin{align}\label{yhf3.16}
\mathbb{E}\left[\exp\bigg(\frac{p}{16}
\sum_{j=1}^{k}\left(8\beta_{r_{j}}+7\lambda\right)
\tr\bigg)\right]\leq H(p)\mathrm{e}^{-\eta_p k \Delta},
\end{align}
and
\begin{align}\label{yhf3.17}
 \begin{split}
 C\tr   \sum_{i=1}^{k} \E\Bigg[ \exp\bigg( \frac{p}{16}\sum_{j=1}^{i}\big(8\beta_{r_{j}}+7\lambda\big)\triangle \bigg) \Bigg]
\leq& C\tr   \sum_{i=1}^{k}H(p)\mathrm{e}^{-\eta_p  i \tr} \\
\leq& C\tr \big(\mathrm{e}^{\eta_p  \tr}-1\big)^{-1}
\leq C.
 \end{split}
\end{align}
Combining  \eqref{yhf3.17}  and  \eqref{yhf3.16}  with  \eqref{111}   yields
\begin{equation}
\begin{split}
   \mathbb  E\Big((1+|X_{k+1}|^{2})^{\frac{p}{2}}\Big)
\leq  C  (1+|x|^{2})^{\frac{p}{2}} \mathrm{e}^{- \eta_p  k \tr}
 +C +C\tr.
 \end{split}
\end{equation}
Therefore the desired assertion follows. 
\end{proof}

\begin{rem}\label{re1}
 {\rm Recently, the work of \cite{Liu15} gives the the moment boundedness of the BEM numerical solutions for  SDEs without globally Lipschitz continuous coefficients. However the proof techniques can't be adopted for SDEs with regime switching directly since their dynamical behaviors are significantly different from those of SDEs. In the proof of Lemma \ref{le-mom} we establish the recursion formula (\ref{eq-tight}) dependent on the states, and then yield the desired result by making use of  the ergodic property of the Markov chain.}
\end{rem}

To investigate  the uniqueness of the invariant measure we provide  the asymptotically attractive property of the numerical solutions of BEM scheme. Here we denote the numerical solution of  BEM scheme with any given initial value $(x,i)$ by  $X_{k}^{x,i}$.
\begin{lemma}\la{le-attr}
Under the conditions of  Theorem \ref{yth3.1}, it holds that
\begin{align}\label{y1}
  \mathbb{E}|X_{k}^{x,i}-X_{k}^{y,j}|^{p}\leq C(1+|x|^p+|y|^p) \mathrm{e}^{-  \varsigma k\triangle}
\end{align}
for any $\tr\in (0, \bar\tr)$ and for any $p\in(0,p_0)$, $(x,i), (y,j)\in \RR^n \times \SS$,
$\bar{\tr}$ and   $p_0$ are given in Lemma \ref{le-mom}, $\varsigma>0$ is a constant.
\end{lemma}
\begin{proof}
Note that
$$ \left\{
\begin{aligned}
X_{k+1}^{x,i}   = & X_{k}^{x,i}+f(X_{k+1}^{x,i},r_{k}^{i})\triangle+g(X_{k}^{x,i},r_{k}^{i})\triangle B_{k}, \\
X_{k+1}^{y,i}  = & X_{k}^{y,i}+f(X_{k+1}^{y,i},r_{k}^{i})\triangle+g(X_{k}^{y,i},r_{k}^{i})\triangle B_{k}.
\end{aligned}\right.
$$
It follows from Assumption \ref{a1} that
\begin{align*}
 |X_{k+1}^{x,i}-X_{k+1}^{y,i}|^{2}
 =&\Big(X_{k+1}^{x,i}-X_{k+1}^{y,i}\Big)^{T}\Big(f(X_{k+1}^{x,i},r_{k}^{i})-f(X_{k+1}^{y,i},r_{k}^{i})\Big)\triangle\\
 &+\Big(X_{k+1}^{x,i}-X_{k+1}^{y,i}\Big)^T\Big(X_{k}^{x,i}-X_{k}^{y,i}+\big(g(X_{k}^{x,i},r_{k}^{i})-g(X_{k}^{y,i},r_{k}^{i})
\big)\triangle B_{k}\Big)\\
  \leq&\alpha_{r_{k}^{i}}\big|X_{k+1}^{x,i}-X_{k+1}^{y,i}\big|^{2}\triangle\\
 &+\frac{1}{2}\big|X_{k+1}^{x,i}-X_{k+1}^{y,i}\big|^{2}+\frac{1}{2}\big|(X_{k}^{x,i}-X_{k}^{y,i})
 +\big(g(X_{k}^{x,i},r_{k}^{i})-g(X_{k}^{y,i},r_{k}^{i})\big)\triangle B_{k} \big|^{2}.
\end{align*}
We hence obtain
\begin{align*}
 |X_{k+1}^{x,i}-X_{k+1}^{y,i}|^{2}
 \leq&\frac{1}{1-2\alpha_{r_{k}^{i}}\triangle}\big|(X_{k}^{x,i}-X_{k}^{y,i})+(g(X_{k}^{x,i},r_{k}^{i})-g(X_{k}^{y,i},r_{k}^{i})\triangle B_{k})\big|^2\\
 =&\frac{|X_{k}^{x,i}-X_{k}^{y,i}|^{2}}{1-2\alpha_{r_{k}^{i}}\triangle}\Big(1+\vartheta(r_{k}^{i})\Big),
\end{align*}
where
\begin{align*}
 \vartheta_{k}(r_{k}^{i})
=&\frac{2(X_{k}^{x,i}-X_{k}^{y,i})^{T}(g(X_{k}^{x,i},r_{k}^{i})-g(X_{k}^{y,i},r_{k}^{i}))\triangle
B_{k}+|\big(g(X_{k}^{x,i},r_{k}^{i})-g(X_{k}^{y,i},r_{k}^{i})\big) \triangle B_{k}|^{2}}{|X_{k}^{x,i}-X_{k}^{y,i}|^{2}}
\end{align*}
if $|X_{k}^{x,i}-X_{k}^{y,i}|\neq0$, otherwise it is set to $-1$.
Clear, $\vartheta_{k}(r_{k}^{i})\geq -1$. For any $p\in (0, p_0)$, then using (\ref{4.5}) we derive that
\begin{align}\label{4.14}
 \begin{split}
  \mathbb E\big(|X_{k+1}^{x,i}-X_{k+1}^{y,i}|^{p}\big|\mathcal F_{t_k}\big)
 \leq&\frac{|X_{k}^{x,i}-X_{k}^{y,i}|^{p}}{\big(1-2\alpha_{r_{k}^{i}}\triangle\big)^{\frac{p}{2}}}I_{\{|X_{k}^{x,i}-X_{k}^{y,i}|\neq0\}}\mathbb E\bigg[1+\frac{p}{2}\vartheta_{k}(r_{k}^{i}) \\
 &~~~~~~~~~~ +\frac{p(p-2)}{8}\vartheta_{k}^{2}(r_{k}^{i})+\frac{p(p-2)(p-4)}{48}\vartheta_{k}^{3}(r_{k}^{i})\big|\mathcal F_{t_k}\bigg].
  \end{split}
\end{align}
Then following the same way as \eqref{4.7}-\eqref{4.9}, by \eqref{a-e2} we can show
\begin{equation}\label{4.15}
  I_{\{|X_{k}^{x,i}-X_{k}^{y,i}|\neq0\}}\mathbb E\big(\vartheta_{k}(r_{k}^{i})|\mathcal F_{t_k}\big)
=  I_{\{|X_{k}^{x,i}-X_{k}^{y,i}|\neq0\}}\frac{|g(X_{k}^{x,i},r_{k}^{i})-g(X_{k}^{y,i},r_{k}^{i})|^{2}\triangle}{|X_{k}^{x,i}-X_{k}^{y,i}|^{2}},
\end{equation}
and
\begin{align}\label{4.16} \begin{split}
  I_{\{|X_{k}^{x,i}-X_{k}^{y,i}|\neq0\}}\mathbb E\big(\vartheta_{k}^{2}(r_{k}^{i})|\mathcal F_{t_k}\big)
 \geq I_{\{|X_{k}^{x,i}-X_{k}^{y,i}|\neq0\}}\frac{4|(X_{k}^{x,i}-X_{k}^{y,i})^{T}(g(X_{k}^{x,i},r_{k}^{i})-g(X_{k}^{y,i},r_{k}^{i}))|^{2}\triangle}{|X_{k}^{x,i}-X_{k}^{y,i}|^{4}},
  \end{split}\end{align}
and
\begin{align}\label{4.17}
  I_{\{|X_{k}^{x,i}-X_{k}^{y,i}|\neq0\}}\mathbb E\big(\vartheta_{k}^{3}(r_{k}^{i})|\mathcal F_{t_k}\big)
 \leq  I_{\{|X_{k}^{x,i}-X_{k}^{y,i}|\neq0\}} C\triangle^{\frac{3}{2}}.
 \end{align}
  Combining (\ref{4.14})-(\ref{4.17}) and
 using   Assumption \ref{a3}, for any $k\geq 0$ we arrive at
\begin{align*}
&\mathbb E\big(|X_{k+1}^{x,i}-X_{k+1}^{y,i}|^{p}|\mathcal F_{t_k}\big)\nn\\
 \leq &\frac{|X_{k}^{x,i}-X_{k}^{y,i}|^{p}}{(1-2\alpha_{r_{k}^{i}}\triangle)^{\frac{p}{2}}}I_{\{|X_{k}^{x,i}-X_{k}^{y,i}|\neq0\}}
 \bigg[1+\frac{p}{2}\bigg(\frac{|g(X_{k}^{x,i},r_{k}^{i})-g(X_{k}^{y,i},r_{k}^{i})|^{2}}{|X_{k}^{x,i}-X_{k}^{y,i}|^{2}}\triangle\\\notag
 &~~~~~~~~~ +(p-2)\frac{|(X_{k}^{x,i}-X_{k}^{y,i})^{T}
 (g(X_{k}^{x,i},r_{k}^{i})-g(X_{k}^{y,i},r_{k}^{i}))|^{2}\triangle}{|X_{k}^{x,i}-X_{k}^{y,i}|^{4}}\bigg)
 +\frac{p(p-2)(p-4)}{48}C\triangle^{\frac{3}{2}}\bigg]\\\notag
 \leq&\frac{|X_{k}^{x,i}-X_{k}^{y,i}|^{p}}{(1-2\alpha_{r_{k}^{i}}\triangle)^{\frac{p}{2}}}
 \bigg[1+\frac{p}{2}\big(h_{r_{k}^{i}}+ph\big)\triangle+\frac{p(p-2)(p-4)}{48}C\triangle^{\frac{3}{2}}\bigg].
\end{align*}
 It is easy to find from (\ref{yp}) that $4ph< \lambda$ holds for each $p\in(0,p_{0})$. Choose a constant $0<\tr_4 \leq \bar{\tr}$ ($\bar{\tr}$ is a positive constant given in Lemma \ref{le-mom}) sufficiently small such that
 $ C\tr_4^{1/2}  \leq {3\lambda }/8, $
 which implies that for any $\tr\in (0,\tr_4]$
\begin{align}\label{eqn6}
\mathbb E\big(|X_{k+1}^{x,i}-X_{k+1}^{y,i}|^{p}|\mathcal F_{t_k}\big)
 \leq \frac{|X_{k}^{x,i}-X_{k}^{y,i}|^{p}}{\big(1-2\alpha_{r_{k}^{i}}\triangle\big)^{\frac{p}{2}}}\bigg[1+\frac{p}{2}\big(h_{r_{k}^{i}}+ \frac{1}{4}\lambda\big)\triangle+\frac{p\lambda }{16} \triangle \bigg].
\end{align}
Further choose  $0<\tr_5 \leq \tr_4$ such that for any $\tr\in (0,\tr_5]$, any $i\in \SS$, any integer $k$
\begin{align}\label{4.18}
(1-2\alpha_{r_{k}^{i}}\triangle)^{\frac{p}{2}}\geq1-p\alpha_{r_{k}^{i}}\triangle-C\triangle^{2}\geq 1-\frac{p}{2}(2\alpha_{r_{k}^{i}}+\frac{1}{16}\lambda)\triangle
\end{align}
 holds. Substituting  this in  \eqref{eqn6} yields
\begin{align*}
\mathbb E\big(|X_{k+1}^{x,i}-X_{k+1}^{y,i}|^{p}|\mathcal F_{t_k}\big)
 \leq \frac{1+\frac{p}{2}( h_{r_{k}^{i}}+\frac{3}{8}\lambda)\triangle}{1-\frac{p}{2}(2\alpha_{r_{k}^{i}}+\frac{1}{16}\lambda)\triangle}|X_{k}^{x,i}-X_{k}^{y,i}|^{p}.
\end{align*}
Using inequality $ {1}/({1-u})\leq 1+u +2u^2$ for any $u\in (-1/2,1/2 )$,  we obtain
\begin{equation}
\begin{split}
 \mathbb E\big(|X_{k+1}^{x,i}-X_{k+1}^{y,i}|^{p}|\mathcal F_{t_k}\big)
\leq \Big(1+\frac{p}{2}(\beta_{r_{k}^{i}}+\frac{1}{2}\lambda)\triangle\Big) |X_{k}^{x,i}-X_{k}^{y,i}|^{p}
\end{split}
\end{equation}
for any $\triangle\in(0,\triangle^*)$, $p\in(0, p_0)$, where $0<\tr^* \leq \tr_5$ satisfying
 $C\triangle^*\leq  {p\lambda}/{32},$ and $ {p_0}(\breve{| \beta|}  + {\lambda}/{2})\triangle^*/2<1. $  This implies that
\begin{align}\label{112}
\begin{split}
\mathbb E\big(|X_{k}^{x,i}-X_{k}^{y,i}|^{p} \big)
 &\leq |x-y|^{p}\mathbb E\bigg[\prod_{j=0}^{k-1}\bigg(1+\frac{p}{2}\Big(\beta_{r_{j}^{i}}+\frac{1}{2}\lambda\Big)\triangle\bigg)\bigg]\\
&\leq|x-y|^{p}\mathbb E\bigg[\exp\bigg ({\sum_{j=0}^{k-1}\log\Big(1+\frac{p}{2}(\beta_{r_{j}^{i}}+\frac{1}{2}\lambda)\triangle\Big)}\bigg)\bigg]\\
&{ \leq|x-y|^{p}\mathbb E\bigg[\exp
\bigg(\frac{p}{4}{\sum_{j=0}^{k-1}\Big(2\beta_{r_{j}^{i}}+\lambda\Big)
\triangle}\bigg)\bigg].}
\end{split}
\end{align}
  For any $p\in (0,p_0)$,
Lemma \ref{yhf2.1}
implies that $\eta_{p}>0$ and there exists  a positive constant  $H(p)$ such that
\begin{align}
\mathbb{E}\left[\exp\bigg(\frac{p}{4}
\sum_{j=0}^{k-1}\left(2\beta_{r^i_{j}}+\lambda\right)
\tr\bigg)\right]\leq\mathbb{E}\left[\exp\bigg(\frac{p}{16}
\sum_{j=0}^{k-1}\left(8\beta_{r^i_{j}}+7\lambda\right)
\tr\bigg)\right]\leq H(p)\mathrm{e}^{-\eta_p k \tr}.
\end{align}
 This  together with (\ref{112})  implies
\begin{align}
\mathbb E|X_{k}^{x,i}-X_{k}^{y,i}|^{p}
\leq H(p) |x-y|^{p}\mathrm{e}^{-\eta_p k \tr},~~~~\forall~k>0.
\end{align}
 Define
$ \bar{\tau} =\inf\{k\geq 0:r_{k}^{i}=r_{k}^{j}\}.$
Since the state space $\SS$ is finite, and $Q$ is irreducible, there exists $\bar{\gamma}>0$ such that
\begin{align}\label{4.21}
\mathbb{P}(\bar{\tau} >k)\leq \mathrm{e}^{-\bar{\gamma} k\triangle}
\end{align}
for any $k>0$. For the fixed $p\in(0,p_0)$, let $q=(p+p_0)/(2p)>1$, then $pq=(p+p_0)/2 \in(0, p_0)$. Moreover,  H\"{o}lder's inequality implies that
\begin{align}\la{h1}
\begin{split}
 &\mathbb E|X_{k}^{x,i}-X_{k}^{y,j}|^{p}\\
 =&\mathbb E\Big(|X_{k}^{x,i}-X_{k}^{y,j}|^{p}I_{\{\bar{\tau} >[\frac{k}{2}]\}}\Big)+\mathbb E\Big(|X_{k}^{x,i}-X_{k}^{y,j}|^{p}I_{\{\bar{\tau} \leq [\frac{k}{2}]\}}\Big)\\
 \leq&\Big(\mathbb E|X_{k}^{x,i}-X_{k}^{y,j}|^{pq}\Big)^{\frac{1}{q}}\Big(\mathbb P(\bar{\tau} >[\frac{k}{2}])\Big)^{1-\frac{1}{q}}+\mathbb E\Big[I_{\{\bar{\tau} \leq [\frac{k}{2}]\}}\mathbb E\Big(|X_{k}^{x,i}-X_{k}^{y,j}|^{p}\big|\mathcal F_{\bar{\tau} \triangle}\Big)\Big]\\
 \leq&\Big(\mathbb E|X_{k}^{x,i}-X_{k}^{y,j}|^{pq}\Big)^{\frac{1}{q}}\Big(\mathbb P(\bar{\tau} >[\frac{k}{2}])\Big)^{1-\frac{1}{q}}+\mathbb E\Big[I_{\{\bar{\tau} \leq [\frac{k}{2}]\}}\mathbb E\Big(|X_{k-\bar{\tau} }^{X_{\bar{\tau} }^{x,i},r_{\bar{\tau} }^i}-X_{k-\bar{\tau} }^{X_{\bar{\tau} }^{y,j},r_{\bar{\tau} }^j}|^{p}\Big)\Big]\\
 \leq&C\mathrm{e}^{-\frac{q-1}{2q}\bar{\gamma} k\triangle}\Big(\mathbb E|X_{k}^{x,i}-X_{k}^{y,j}|^{pq}\Big)^{\frac{1}{q}}
 +C{\mathrm{e}^{-\frac{\eta_p}{2}  k\triangle}}\mathbb E\Big[I_{\{\bar{\tau} \leq [\frac{k}{2}]\}}\mathbb E\Big(|X_{\bar{\tau} }^{x,i}-X_{\bar{\tau} }^{y,j}|^{p}\Big)\Big]\\
 \leq& C \mathrm{e}^{-\frac{p_0-p}{2(p +p_0)}\bar{\gamma} k\triangle}\Big(\mathbb E|X_{k}^{x,i}-X_{k}^{y,j}|^{\frac{p+p_0}{2}}\Big)^{\frac{2p}{p+p_0}}+C\mathrm{e}^{-\frac{\eta_p}{2}  k\triangle}\mathbb E\Big(\big|X_{\bar{\tau} \wedge [\frac{k}{2}]}^{x,i}-X_{\bar{\tau} \wedge [\frac{k}{2}]}^{y,j}\big|^{p}\Big),
 \end{split}
\end{align}
where $[x]$ represents the integer part of $x$ for any $x\in \RR$.
  Applying the elementary inequality $(a+b)^p\leq 2^p(a^p+b^p)$ for all $a, b>0$, {by (\ref{4.13})},  yields that
$$\big(\mathbb E|X_{k}^{x,i}-X_{k}^{y,j}|^{{\frac{p+p_0}{2}}}\big)^{\frac{2p}{p+p_0}}\leq C(1+|x|^p+|y|^p),$$ and
\begin{align*}
 \mathbb E\Big(|X_{\bar{\tau} \wedge [\frac{k}{2}]}^{x,i}-X_{\bar{\tau} \wedge [\frac{k}{2}]}^{y,j}|^{p}\Big)
\leq& \mathbb E\Big(|X_{\bar{\tau} \wedge [\frac{k}{2}]}^{x,i}|^{p}\Big)+\mathbb E\Big(|X_{\bar{\tau} \wedge [\frac{k}{2}]}^{y,j}|^{p}\Big)\nn\\
=&   \mathbb E\Big(\sum_{l=0}^{[\frac{k}{2}]}|X_{l}^{x,i}|^{p}I_{\{\bar{\tau} \wedge [\frac{k}{2}]=l\}}(\omega)\Big)+ \mathbb E\Big(\sum_{l=0}^{[\frac{k}{2}]}|X_{l}^{y,j}|^{p}I_{\{\bar{\tau} \wedge [\frac{k}{2}]=l\}}(\omega)\Big)\nn\\
\leq& \sum_{l=0}^{[\frac{k}{2}]}\Big[\mathbb E\big(|X_{l}^{x,i}|^{p}\big)+ \mathbb E\big(|X_{l}^{y,j}|^{p}\big)\Big]
\leq     C(1+|x|^p+|y|^p)(k+2 ).
\end{align*}
  The desired assertion (\ref{y1}) follows by using  (\ref{h1}).
  \end{proof}

Next we give the existence and uniqueness of the numerical invariant measure for SDE (\ref{1.1}) of  BEM method.
\begin{theorem}\la{yth3.2}
Under the conditions of  Theorem \ref{yth3.1},
 there is a positive $\triangle^{*} $
 sufficiently small such that for any $\triangle\in(0,\triangle^{*})$,
 the solutions of the BEM method (\ref{4.2}) converge to a unique invariant measure $\pi^{\Delta}\in \mathcal{P}(\RR^n \times \SS)$ with some exponential rate $\xi_\tr>0$ in the Wasserstein distance.
\end{theorem}
\begin{proof}
 For any initial data $(x,i)$, by (\ref{4.13}) and Chebyshev's inequality, we derive that $\{\delta_{(x,i)}\mathbf{P}_{k\triangle}^{\triangle}\}$ is tight, then one can extract a subsequence which converges weakly to an invariant measure denoted by $\pi^{\Delta}\in \mathcal{P}(\RR^n \times \SS)$.
It follows from  (\ref{4.21}) that
\begin{align}\la{y2}
\mathbb{P}(r_{k}^{i}\neq r_{k}^{j})=\mathbb{P}(\bar{\tau} >k)\leq \mathrm{e}^{-\bar{\gamma} k\triangle}
\end{align}
for any $k>0$.   Therefore, we derive  from (\ref{y1}) and (\ref{y2})  that
\begin{align}\la{y2.15}\begin{split}
W_{p}(\delta_{(x,i)}\mathbf{P}_{k\triangle}^{\triangle},\delta_{(y,j)}
\mathbf{P}_{k\triangle}^{\triangle})\leq& \mathbb{E}|X_{k}^{x,i}-X_{k}^{y,j}|^{p}+\mathbb{P}(r_{k}^{i}\neq r_{k}^{j})\\ \leq& C(1+|x|^p+|y|^p) \mathrm{e}^{-  \xi_\tr k\triangle}.
\end{split}
\end{align} where $\xi_\tr:= \varsigma \wedge\bar{\gamma}>0.$
Due to the Kolmogorov-Chapman equation and Lemma \ref{le-mom} one observes that for any $k,~l>0$,
\begin{align*}
 W_{p}(\delta_{(x,i)}\mathbf{P}_{k\triangle}^{\triangle},\delta_{(x,i)}
\mathbf{P}_{(k+l)\triangle}^{\triangle})
=& W_{p}(\delta_{(x,i)}\mathbf{P}_{k\triangle}^{\triangle},\delta_{(x,i)}
\mathbf{P}_{l\triangle}^{\triangle}\mathbf{P}_{k\triangle}^{\triangle})\\
 \leq &\int_{\RR^{n}\times \SS} W_{p}(\delta_{(x,i)}\mathbf{P}_{k\triangle}^{\triangle},\delta_{(y,j)}
 \mathbf{P}^{\triangle}_{k\triangle})\mathbf{P}_{l\triangle}^{\triangle}(x,i; dy, j)\\
\leq &\sum_{j\in\SS}\int_{\RR^{n} }   C(1+|x|^p+|y|^p) \mathrm{e}^{-   \xi_\tr k\triangle}\mathbf{P}_{l\triangle}^{\triangle} (x,i;dy,j)\\
=& C(1+|x|^p+\mathbb{E}|X_{l}^{x,i} |^{p}) \mathrm{e}^{-   \xi_\tr k\triangle} \leq C  \mathrm{e}^{-  \xi_\tr k\triangle}.
\end{align*} Thus, taking $l\rightarrow\infty$ implies
\begin{equation}\label{eq-expo}
W_{p}(\delta_{(x,i)}\mathbf{P}_{k\triangle}^{\triangle}, \pi^{\Delta})\leq C  \mathrm{e}^{-  \xi_\tr k\triangle}\rightarrow 0,~~~~k\rightarrow\infty,
\end{equation}
namely, $\pi^{\Delta}$ is the unique invariant measure of $\{\delta_{(x,i)}\mathbf{P}_{k\triangle}^{\triangle}\}$.
Assume   $\nu_1^{\triangle}, \nu_2^{\triangle}\in \mathcal{P}(\RR^n \times \SS)$ are the invariant measures of $(X_k^{x,i},r_k^i)$ and $(X_k^{y,j},r_k^j)$, respectively,
 we have
\begin{align*}
W_{p}(\nu_1^{\triangle},\nu_2^{\triangle}) \le \int_{(\RR^{n}\times \SS) \times(\RR^{n}\times \SS)} W_{p}(\delta_{(x,i)}\mathbf{P}_{k\triangle}^{\triangle},\delta_{(y,j)}
 \mathbf{P}^{\triangle}_{k\triangle})\pi(dx\times di,dy\times dj),
\end{align*}
where $\pi$ is a coupling of $\nu_1^{\triangle}$ and $\nu_2^{\triangle} $.
Therefore, the uniqueness of invariant measures follows from (\ref{y2.15})  immediately. \end{proof}

The following theorem reveals that numerical invariant measure $\pi^{\triangle}$ converges in the Wassertein distance to the underlying one $\pi$.
\begin{theorem}\la{yth3.3}
 Under the conditions of  Theorem \ref{yth3.1},
$
\lim_{\triangle\to 0}W_{p}(\pi,\pi^{\triangle})=0.
$ Furthermore, if the drift term satisfies the  polynomial growth condition, that is,
$$|f(x,i) -f(y,i)|^2\leq c_i(1+|x|^q+|y|^q)|x-y|^2,  ~~\forall x,y\in \RR^n, i\in \SS,$$
then $
 W_{p}(\pi,\pi^{\triangle})\leq  C \tr^\gamma
$ for some $\gamma\in (0,p/2)$, where $c_i,~q$ are positive constants.
\end{theorem}
\begin{proof}  Under Assumptions \ref{a1} and \ref{a3}, by Theorem \ref{yth3.1}, Remark \ref{re3} and  \eqref{eq-expo}, for any $\triangle\in (0, \tr^*)$ and any $\epsilon>0$, there is a $k>0$ sufficiently large  such that
\begin{align}\label{cf1}
W_{p}(\delta_{(x,i)}\mathbf{P}_{k\triangle},\pi)
+W_{p}(\delta_{(x,i)}\mathbf{P}_{k\triangle}^{\triangle},\pi^{\triangle})
\le C  \mathrm{e}^{-  \xi^* k\triangle }<\frac{\epsilon}{2},
\end{align}
 where $\tr^*$ is given by Theorem \ref{yth3.2} and  $\xi^*:=\xi  \wedge \xi_\tr$. Moreover, for the fixed $k$ by the convergence of finite time  when $\triangle$ is sufficiently small,
$$
 W_{p}(\delta_{(x,i)}\mathbf{P}_{k\triangle},
 \delta_{(x,i)}\mathbf{P}_{k\triangle}^{\triangle})<\frac{\epsilon}{2}.
$$
Therefore
the first desired assertion follows.

Furthermore, under the polynomial growth condition of $f$, by the similar way to \cite{H01}, we can obtain that
$$ W_{p}(\delta_{(x,i)}\mathbf{P}_{k\triangle},\delta_{(x,i)}\mathbf{P}_{k\triangle}^{\triangle})\leq Ce^{\nu k\triangle}\tr^{p/2},$$
for some positive constant $\nu$.
 Let $\bar{K}$ be the integer part of  constant $ - {p\ln\tr}/[{2(\nu+\xi^*)\tr} ]$, obviously, $\bar{K}\rightarrow \infty$ as $\tr\rightarrow 0$. One observes that
$$ e^{\nu \bar{K}\triangle}\tr^{p/2}\leq \tr^{\frac{ p \xi}{2(\nu+\xi^*)}},~~~\mathrm{e}^{-  \xi \bar{K}\triangle }\leq   e^{\xi^* \tr^*}\tr^{\frac{  p \xi}{2(\nu+\xi^*)}}.$$
 Therefore,
 $W_{p}(\pi, \pi^{\triangle})\leq  C \tr^{\frac{p \xi }{2(\nu+\xi^*)}}=:C \tr^\gamma.
$ \end{proof}

\begin{rem}\la{re2}
 {\rm In Theorem \ref{yth3.3} we  not only give the convergence of invariant measures but also reveal the rate of the convergence is  exponential under the polynomial growth condition imposed on  $f$. We also notice that Meyn and Tweedie's work  \cite{Me1} reveals the relationship of tightness, Harris recurrence and ergodicity for discrete-time Markov chains, they gave the generalization of Lyapunov-Foster criteria for the various ergodicity. However, these criteria are not applicable for $(X_{k},r_{k})$ owing to the switching effects.  Precisely,  it is impossible  from (\ref{eq-tight}) to find a constant $0<\lambda\leq 1$ such that
$
 \mathbb E\big((1+|X_{k+1}|^{2})^{\frac{p}{2}}|\mathcal F_{t_k}\big)
  \leq \lambda(1+|X_{k}|^{2})^{\frac{p}{2}}+C\triangle$ holds due to the changeable sign of $
  \beta_{r_{k}}+ {7} \lambda/8$.}
   \end{rem}

   \begin{rem} {\rm
  By the virtue of Theorem \ref{yth3.2}, $(X_{k},r_{k})$ is ergodic, moreover,
    the transition probability of  $(X_{k},r_{k})$ decays into its invariant measure exponentially under Wasserstein distance, see \eqref{eq-expo}.}
 \end{rem}

\begin{rem}\la{rem4}
{\rm Comparing with the  convergence result of the EM scheme for SDE in \cite{Shard+stu}, we release the restriction of  the global Lipschitz continuity of the coefficients and deal with the convergence of invariant measures for nonlinear SDE with regime switching.} \end{rem}

\begin{rem}
{\rm  Although many works pay attention to the approximation of invariant measures for SDEs, for example, \cite{Liu15, Ma02, Shard+stu}, there are few works focusing on the approximation of invariant measures for  switching diffusion processes, especially described by nonlinear systems. On the other hand, compared with the fast development of the finite-time numerical analysis for SPDEs, for examples, \cite{Je, Shard}, the results on long-time approximations for SPDEs are few. The methods developed in this paper provide  ideas to  deal with the invariant measure approximations for nonlinear SPDEs or  SPDEs  with regime switching. Owing to the importance this will be considered in our future work.}
 \end{rem}

\section{Examples}\label{sec:example}
In this section, we consider two examples of nonlinear hybrid stochastic systems and provide  simulations to illustrate the efficiency of the  BEM method (\ref{4.2}).    We first consider a two-dimensional  SDE with Markovian switching.

\begin{expl}\la{exp1}{\rm
Consider (\ref{1.1}) with $r(t)$ taking values in $\SS=\{1, 2 \}$ with generator $$Q=\left(
  \begin{array}{ccc}
    -5 & 5\\
    1 & -1\\
  \end{array}
\right).
$$  The  system is regarded  as the Markovian switching between
\begin{align}\la{BEM1}
\left\{
\begin{array}{ll}
\mathrm{d}Y_1(t)=\big[2Y_1(t)-Y_1^3(t)-Y_1(t)Y_2^2(t)\big] \mathrm{d}t-3\mathrm{d}B_1(t)+\mathrm{d}B_2(t),\\
\mathrm{d}Y_2(t)=\big[1+Y_2(t)-Y_2^3(t)-Y_2(t)Y_1^2(t)\big] \mathrm{d}t+4\mathrm{d}B_1(t),
\end{array}
\right.
\end{align}
and
\begin{align}\la{BEM2}
 \left\{
\begin{array}{ll}
\mathrm{d}Y_1(t) &= \Big(Y_1(t)\!-2Y_1(t)\sqrt{Y_1^2(t) +Y_2^2(t)}+1\Big) \mathrm{d}t\\
 &~~~~+(2Y_1(t) -Y_2(t)+2)\mathrm{d}B_1(t) +(Y_1(t) -Y_2(t))\mathrm{d}B_2(t),\\
\mathrm{d}Y_2(t)&=\Big(0.5Y_2(t)-2Y_2(t)\sqrt{Y_1^2(t) +Y_2^2(t)}+2\Big) \mathrm{d}t\\
 &~~~~+(Y_1(t) +2Y_2(t))\mathrm{d}B_1(t) +(Y_1(t)+Y_2(t)-4)\mathrm{d}B_2(t),
\end{array}
\right.
\end{align}
with  the initial data $Y(0)=1$, $r(0)=1$, where $B(t)=(B_1(t),B_2(t))^T$  is a two-dimensional Brownian motion.  Obviously, the diffusion coefficient  $g$ is global Lipschitz continuous with $h=7$.
  Note that the drift coefficient
$
f
$
  is neither the global Lipschitz continuous nor the linear growth, but
 we can derive that
 $$
(u-v)^{T}(f(u,1)-f(v,1))\leq 2|u-v|^{2}-\frac{1}{4}(|u|-|v|)^4\leq 2|u-v|^{2},
$$
and
$$
(u-v)^{T}(f(u,2)-f(v,2))\leq  |u-v|^{2}-2(|u|+|v|)(|u|-|v|)^2\leq  |u-v|^{2},
$$
i.e. Assumption \ref{a1} is satisfied with $\alpha_1=2$ and $\alpha_2=1$ for all $u, v\in \mathbb{R}^2$.  We furthermore
 observe  that 
\begin{align*}
&|u-v|^{2}|g(u,j)-g(v,j)|^{2}-2|(u-v)^{T}(g(u,j)-g(v,j))|^{2}\leq h_j |u-v|^{4},~~~~\forall j\in \mathbb{S},
\end{align*}
holds with $h_1=0$ and $h_2=-3$ for all $u, v\in \mathbb{R}^2$.  Direct calculation leads to
$
\beta_1=2\alpha_1+h_1=4,~ \beta_2=2\alpha_2+h_2=-1.
$
 By solving the linear equation (\ref{eq:a1.2}) we obtain the unique stationary distribution of $r(t)$,
$
\mu=\left(\mu_1, \mu_2\right)=\left( {1}/{6},  {5}/{6}\right),
$  then $
\mu \beta=\mu_1 \beta_1+\mu_2 \beta_2=- {1}/{6}<0.
$
It follows  from Theorem \ref{yth3.1} that the  exact solution $(Y(t),r(t))$ of (\ref{1.1}) admits a unique invariant measure $\pi\in \mathcal{P}(\RR^n\times \SS)$.
By virtue of Theorems \ref{yth3.2} and \ref{yth3.3}, for a given stepsize $\triangle$ the numerical solution of BEM scheme has a unique invariant measure $\pi^{\tr}\in \mathcal{P}(\RR^n\times \SS)$ approximating $\pi$ in the Wasserstein metric. We apply  the BEM    scheme for numerical experiments. Since it is impossible to get the closed form of the solutions of  the stochastic system with random switching  between (\ref{BEM1}) and (\ref{BEM2}), we approximate the underlying solution by the numerical solution of BEM   scheme (\ref{4.2}). We regard the numerical solution with $\triangle=2^{-17}$ as a more precise  approximation comparing it with the numerical solution with  stepsize $\triangle=0.002$, see Figure \ref{figure1}.
\begin{figure}[!htbp]
  \centering
\includegraphics[width=16cm,height=9cm]{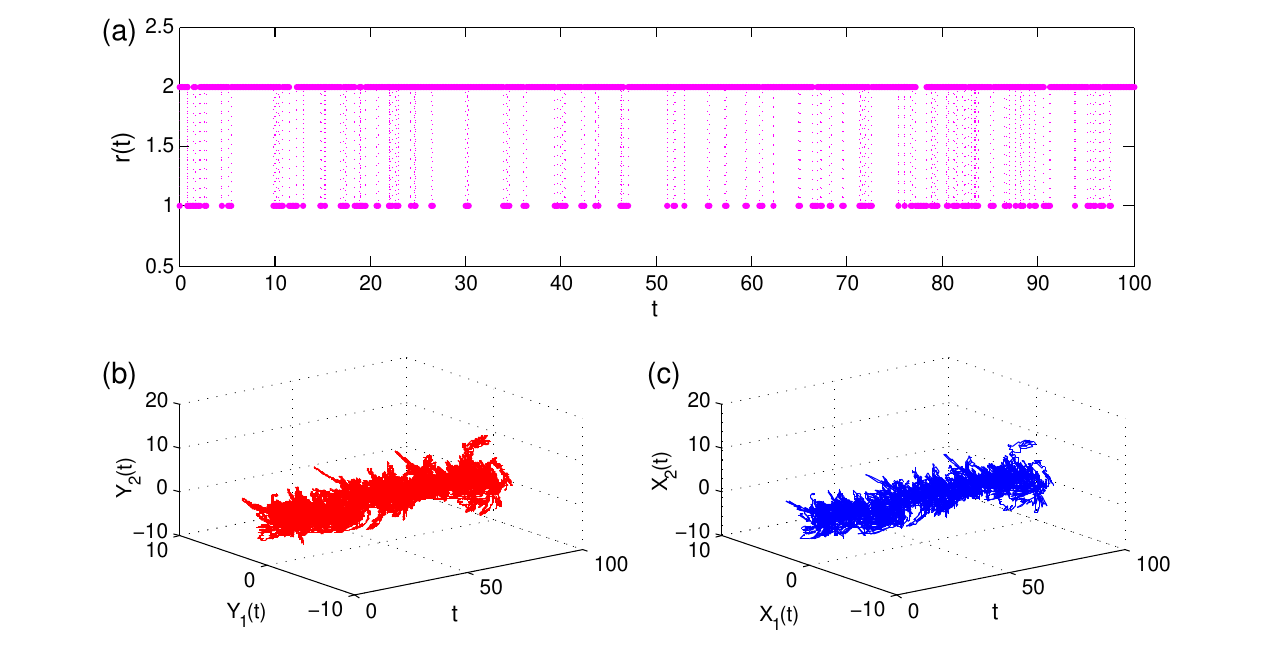}
  \caption{$\mathbf{Example~\ref{exp1}}$. (a)  Computer simulation of a single path of  Markov chain $r(t)$.  (b) A sample path of exact solution $Y(t)$ in 3D settings. (c) A sample path of numerical solution $X(t)$ in 3D settings. The red  line represents the exact solution (i.e. the BEM numerical solution  with $\triangle=2^{-17}$) while  the blue line represents the BEM numerical solution with $\triangle=0.002$.}
  \label{figure1}
\end{figure}
\begin{figure}[!htbp]
  \centering
\includegraphics[width=16cm,height=6cm]{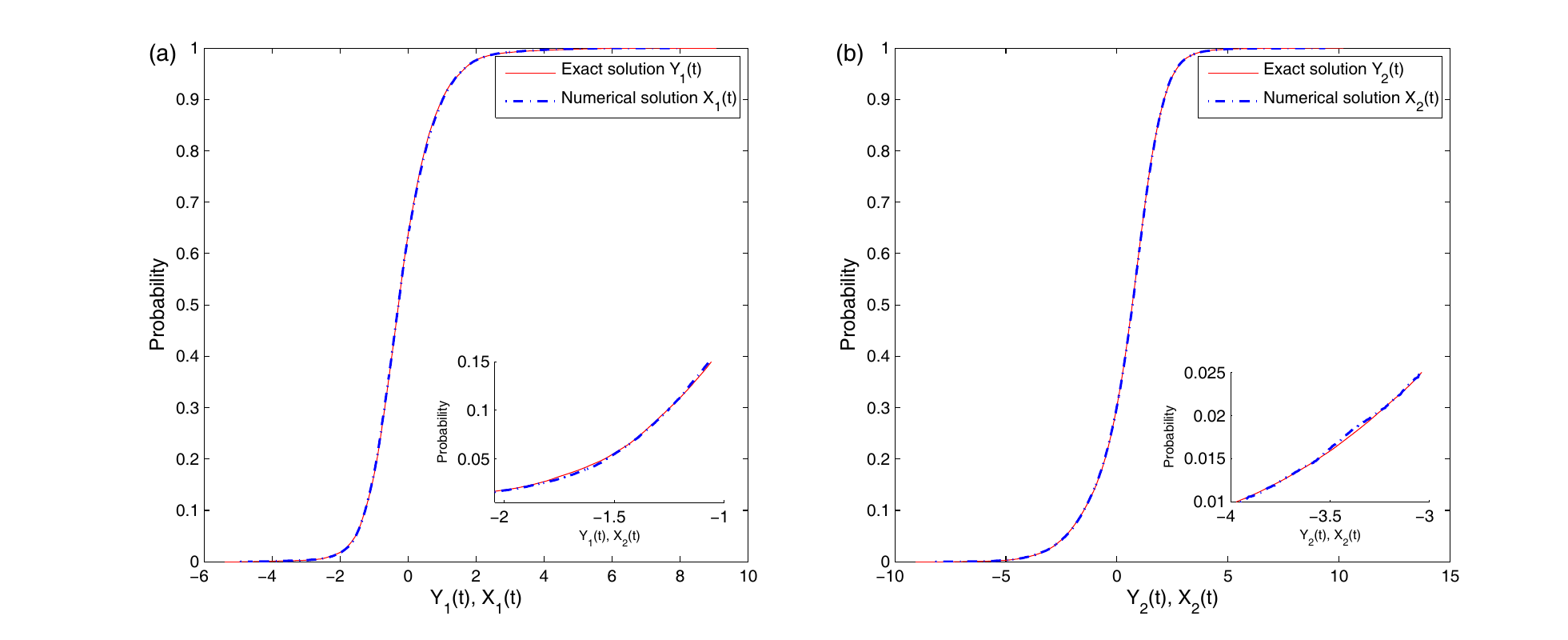}
  \caption{$\mathbf{Example~\ref{exp1}}$. (a)  The ECDF for $Y_1(t)$.  (b) The ECDF for $Y_2(t)$.   The red solid line represents the exact solution of the switching system while the blue dashed line represents the numerical solution of the switching system.}
  \label{figure2}
\end{figure}

We simulate one path with 13107200 iterations and plot the empirical cumulative distribution function (ECDF) of numerical solution with $\triangle=0.002$ in blue dashed line in Figure \ref{figure2}. The ECDF of exact solution is plotted on the same figure in a red solid line.
 The similarity of those two distributions is clearly seen, which indicates that the numerical stationary distribution is a good approximation to the theoretical one.
  To measure the similarity quantitatively, we use the Kolmogorov-Smirnov test  \cite{Massey} to test the alternative hypothesis that the numerical solution and exact solution are from different distributions against the null hypothesis that they are from the same distribution for both $Y_1(t)$ and $Y_2(t)$. With 3\% significance level, the  Kolmogorov-Smirnov test
indicates that we cannot reject the null hypothesis. This example illustrates that numerical invariant measure  converges to the underlying invariant measure.

}
\end{expl}

In order to illustrate the validity,  we consider   the scalar hybrid cubic  SDE
 (c.f.  the stochastic Ginzburg-Laudau equation (4.52) in \cite[p.125]{Kloeden} ) which drift coefficient isn't global Lipschitz continuous. 

\begin{expl}\la{exp2}{\rm  Let $r(t)$ be a Markov chain with the state space $\SS=\{1, 2 \}$ and the generator  $$Q=\left(
  \begin{array}{ccc}
    -q & q\\
    3 & -3\\
  \end{array}
\right)
$$  for some $q>0$.  It is easy to see that its unique stationary distribution $\mu=\left(\mu_1, \mu_2\right)\in \mathbb{R}^{1\times 2}$ is given by
$
 \mu_1=\frac{3}{3+q}, ~ \mu_2=\frac{q}{3+q}.
$
 Consider the scalar hybrid cubic SDE
\begin{align}\la{exyhf4.1}
\mathrm{d}Y(t)=(b(r(t))Y(t)+a(r(t)) Y^3(t) )\mathrm{d}t+\rho(r(t)) Y(t)\mathrm{d}B(t),
\end{align}
 with the initial data $Y(0)=0.5$, $r(0)=2$, where
 $$
b(1)=1,~~a(1)=-1,~~\rho(1)=2,~~~~~b(2)=2,~~a(2)=-3,~~\rho(2)=-1,
 $$
and $B(t)$ is a scalar Brownian motion. There exists a unique continuous  solution $Y(t)$ to SDE \eqref{exyhf4.1} for any $Y(0)>0$, which is global and represented by
\begin{align*}
 Y(t)= \frac{0.5\exp{\Big\{\dis\int_{0}^{t} \Big[b(r(s))-\frac{1}{2}\rho^2(r(s))\Big]\mathrm{d}s+\rho(r(s))\mathrm{d}B(s)\Big\}}}
 {\sqrt{1-0.5\dis\int_{0}^{t}a(r(s))\exp{\Big\{\int_{0}^{s} \Big[2b(r(u))- \rho^2(r(u))\Big]\mathrm{d}u+2\rho(r(u))\mathrm{d}B(u)\Big\}}}\mathrm{d}s}.
\end{align*}
  It is straightforward to see that $\alpha_1=1$, $\alpha_2=2$, $h_1=-4$, and $h_2=-1$.    Direct calculation leads to
$
\beta_1=-2,~ \beta_2=3,
$
 then $$
\mu \beta=\mu_1 \beta_1+\mu_2 \beta_2<0
$$
holds with $q\in(0, 2)$. It follows  from Theorem \ref{yth3.1} that the  exact solution $(Y(t),r(t))$ of (\ref{1.1}) admits a unique invariant measure $\pi\in \mathcal{P}(\RR^n\times \SS)$.
By virtue of Theorems \ref{yth3.2} and \ref{yth3.3} the numerical solution of BEM scheme has a unique invariant measure $\pi^{\tr}\in \mathcal{P}(\RR^n\times \SS)$ approximating $\pi$ in the Wasserstein metric.
\begin{figure}[!htbp]
  \centering
\includegraphics[width=16cm,height=6.8cm]{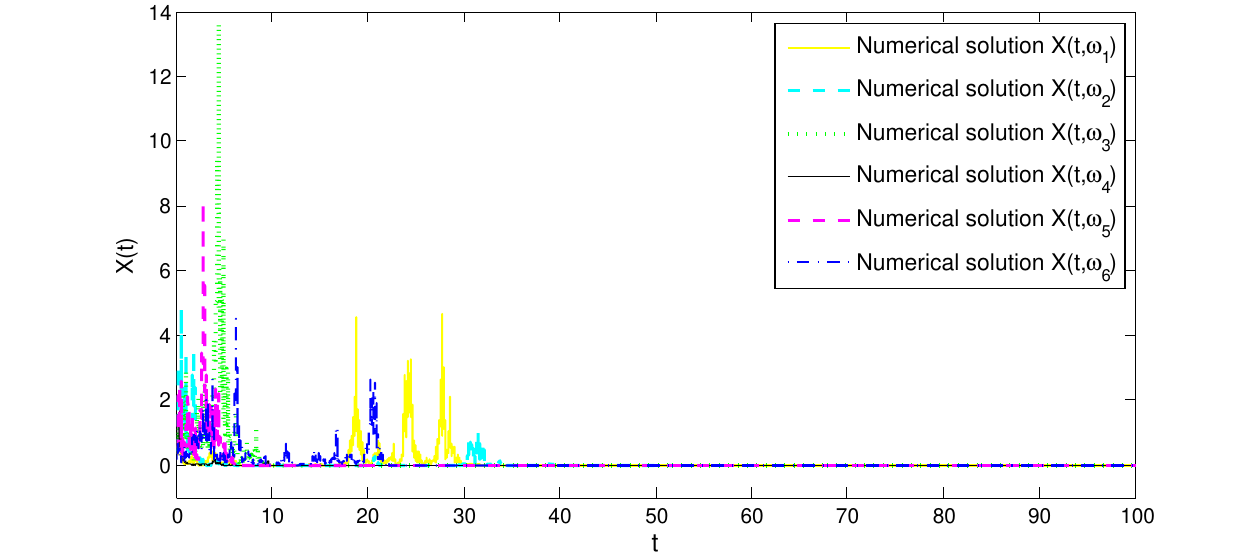}
  \caption{$\mathbf{Example~\ref{exp2}}$. Six trajectories of the BEM numerical solution with $10^4$ iterations, $Y(0)=0.5$, $r(0)=2$ and stepsize  $\triangle=0.001$.}
  \label{figure3}
\end{figure}
\begin{figure}[!htbp]
  \centering
\includegraphics[width=16cm,height=8cm]{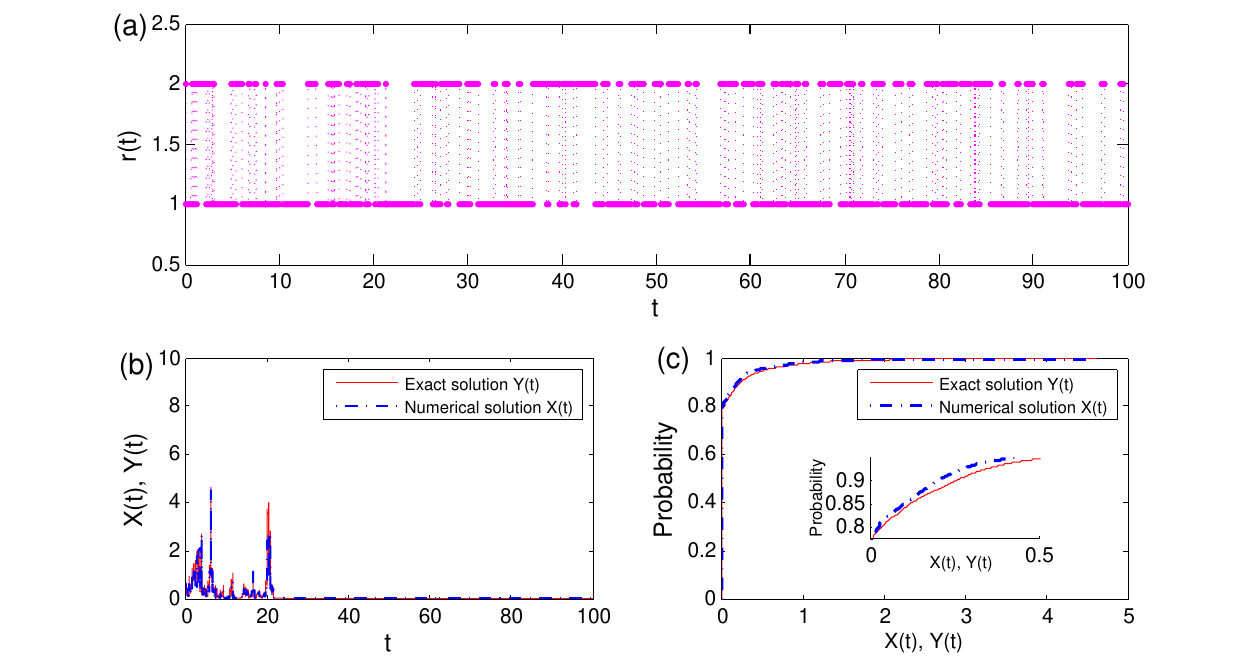}
  \caption{$\mathbf{Example~\ref{exp2}}$. (a)  Computer simulation of a single path of  Markov chain $r(t)$.  (b) Sample paths of the exact solution and the BEM solution. (c) ECDFs for the exact solution and the BEM solution. The red solid line represents exact solution of the switching system while the blue dashed line represents the  numerical solution.}
  \label{figure4}
\end{figure}

We apply  the BEM  scheme to do numerical experiments.
 Choose $q=1.5$ and  stepsize  $\triangle=0.001$,  we simulate $100$ paths, each of which has $10^4$ iterations. Figure \ref{figure3} depicts 
 six trajectories of the numerical solution of BEM   scheme (\ref{4.2}). Intuitively, some stationary behaviours display.
Figure \ref{figure4} (a) depicts the trajectory of the Markov chain. From this figure we find that the time  the Markov chain staying on state 1   is more than on that of state 2.
Figure \ref{figure4}(b) further depicts the trajectories of the exact solution $Y(t)$ and the corresponding BEM solution $X(t)$, and    Figure \ref{figure4}(c) depicts the ECDFs of the  exact  solution and the  BEM   solution.  The similarity of those two distributions is clear, which reveals that the numerical stationary distribution is a good approximation to the underlying one.  Moreover,
 This example illustrates the existence of the stationary distribution as
time goes to infinity. Thus instead of using numerous paths, we could just use few paths to picture the stationary distribution.
}
\end{expl}

\section*{Acknowledgements}

The authors would like to thank the associated editor and referee for their
helpful comments and suggestions.

\bibliography{siam_arxiv}
\end{document}